\documentclass[12pt]{article}

\usepackage{amsthm,amsfonts,latexsym,amssymb}

\theoremstyle{plain}
\newtheorem{proposition}{Proposition}
\newtheorem*{theorem*}{Theorem}

\theoremstyle{definition}
\newtheorem*{definition}{Definition}

\def\pf{\ifvmode\else\newline\fi\noindent\textsc{Proof:\ }}
\def\qed{\mbox{ $\Box$}}

\def\D{{\mathbb D}}
\def\S{{\mathbb S}}

\def\txt#1{\quad\mbox{#1}\quad}

\def\tsprod#1#2{\langle #1, #2\rangle}
\def\eqref#1{\mbox{\rm (\ref{#1})}}
\def\UMD{UMD}
\def\colon{\;:\;}
\def\m{\mbox{\boldmath$\mu$}}
\def\tfrac#1#2{{\textstyle \frac{#1}{#2}}}
\def\spec{\mathop{\rm spec}}

\renewenvironment{cases}{\left\{\begin{array}{cl}}{\end{array}\right.}

\begin{document}
\def\footnotemark{}
\date{April 25, 1997}
\title{Ideal norms associated with the \UMD-property\thanks{Research
    supported by German Academic Exchange Service (DAAD)}}
\author{J\"org Wenzel}
\maketitle

\begin{abstract}
  We prove the asymptotic equivalence of three sequences of ideal
  norms associated with the \UMD-property of Banach spaces.
\end{abstract}

\section{Introduction}
\label{sec:intro}

A Banach space $X$ is a \UMD-space, if there is a constant $c\geq 1$
such that
\begin{equation} \label{eq:umd}
  \Bigg(
    \int_M \bigg\| \sum_{k=1}^n \epsilon_k d_k(\xi) \bigg\|^2 d\mu(\xi)
  \Bigg)^{1/2} \leq c
  \Bigg(
    \int_M \bigg\| \sum_{k=1}^n d_k(\xi) \bigg\|^2 d\mu(\xi)
  \Bigg)^{1/2}
\end{equation}
for all sequences $d_1,\dots,d_n$ of $X$-valued martingale differences
and all sequences $\epsilon_1,\dots,\epsilon_n$ of signs.  (The
letters \UMD{} stand for \emph{unconditional martingale differences}.)
Maurey \cite{mau74} and later Burkholder \cite{bur86} showed, that
this is the case if and only if \eqref{eq:umd} is satisfied for
Walsh-Paley-martingales on the interval $[0,1)$ only.  Throughout this
article, we will only deal with those special martingales.

In this setting, there are essentially three different ways of
changing signs:
\begin{enumerate}
\item use all \emph{predictable} sequences $(\epsilon_k)$, i.~e.
  $\epsilon_k \colon [0,1]\to\{\pm1\}$ is ${\cal F}_{k-1}$-measurable,
  where $({\cal F}_k)$ is the filtration, to which the martingale is
  adapted,
\item use all constant sequences of signs $\epsilon_k\in\{\pm1\}$,
\item use one fixed sequence of signs $\epsilon_k=(-1)^k$.
\end{enumerate}

For each fixed $n$ in \eqref{eq:umd}, we will define below three
corresponding ideal norms. The obtained sequences of ideal norms are
bounded, if and only if $X$ is a \UMD-space. However, also in the
non-bounded case we can gain some information on $X$ from the
asymptotic behavior of these sequences. The main result of this paper
states that this information is essentially the same in all three
cases. The corresponding sequences of ideal norms are asymptotically
equivalent.

A similar result in the setting of general martingales was obtained by
Burkholder in \cite[Lemma 2.1]{bur84}. However, to make his proof
work, one has to allow the underlying filtrations for the martingales
to vary.

In the natural way, all concepts extend to the setting of operators
between Banach spaces.

\section{Definitions and main result}
\label{sec:def}

For $k=1,2,\dots$ and $j=0,\pm1,\pm2,\dots$, we let
\[
  \Delta_k^{(j)}:=\left[\tfrac{j-1}{2^k},\tfrac j{2^k}\right)
\]
be the {\em dyadic intervals}. The {\em Haar functions} are given by
\[
  \chi_k^{(j)}(t):=
  \begin{cases}
    +2^{(k-1)/2} & \txt{if $t\in\Delta_k^{(2j-1)}$,}\\[3pt]
    -2^{(k-1)/2} & \txt{if $t\in\Delta_k^{(2j)}$,}\\
    0 & \txt{otherwise.}
  \end{cases}
\]

We let
\[
  \D:=\{(k,j)\colon k=1,2,\dots;\ j=1,\dots,2^{k-1}\}
\]
denote the {\em dyadic tree}. We will mainly consider finite dyadic
trees
\[
  \D_m^n:=\{(k,j)\colon k=m,\dots,n;\ j=1,\dots,2^{k-1}\},
\]
where $m\le n$. To shorten terms, we write $\D_k$ for the \emph{$k$-th
  level} $\D_k^k$ of $\D$.

We denote by $L_2^X$ the Banach space of square integrable $X$-valued
functions $f$ on the interval $[0,1)$ equipped with the norm
\[
\|f\|_2 := \bigg( \int_0^1 \|f(t)\|^2 dt \bigg)^{1/2}.
\]
All results in this article could also be obtained for an arbitrary
index $1 < p < \infty$ instead of $2$, the changes are
straightforward. However, to avoid cumbersome notation, we decided to
restrict ourselves to the case $p=2$.

Given any $\D_1^n$-tuple $(x_k^{(j)})$, we get a
\emph{Walsh-Paley-Martingale} of length $n$ with mean value zero, by
letting
\[
f_k := \sum_{(h,i)\in\D_1^k} \!\! x_h^{(i)} \chi_h^{(i)}
\txt{for $k=1,\dots,n$.}
\]
Note that by the martingale properties of the sequence $(f_k)$ and
since the conditional expectation operator has norm one in $L_2^X$, we
have
\begin{equation}
  \label{eq:cond_exp}
  \|f_k\|_2 \leq \|f_n\|_2
\end{equation}
whenever $k\leq n$.  We write
\[
\tsprod f{\chi_k^{(j)}} := \int_0^1 f(t) \chi_k^{(j)}(t) \, dt,
\]
for the \emph{Haar-Fourier coefficients} of a function $f\in L_2^X$
and call
\[
\spec(f) := \{ (k,j)\in\D \colon \tsprod f{\chi_k^{(j)}} \not=0\}
\]
the \emph{spectrum} of the function $f$.

\begin{definition}
  For an operator $T\colon X \to Y$, we denote by $\m_n(T)$ the least
  constant $c\geq 1$ such that
  \[
  \bigg\|
  \sum_{(k,j)\in\D_1^n} \!\! \epsilon_k^{(j)} Tx_k^{(j)} \chi_k^{(j)}
  \bigg\|_2
  \leq c \,
  \bigg\|
  \sum_{(k,j)\in\D_1^n} \!\! x_k^{(j)} \chi_k^{(j)}
  \bigg\|_2
  \]
  for all $\D_1^n$-tuples $(x_k^{(j)})$ and all signs
  $\epsilon_k^{(j)}=\pm1$.
\end{definition}

The above definition can be modified by assuming that the signs are
changed on every level simultaneously. In other terms,
$\epsilon_k^{(j)}=\epsilon_k=\pm 1$ should not depend on
$j=1,\dots,2^{k-1}$. A still weaker concept can be introduced by using
only the signs $\epsilon_k^{(j)}=(-1)^k$.  The ideal norms so obtained
will be denoted by $\m_n^\circ(T)$ and $\m_n^{\circ\circ}(T)$,
respectively.

Note that the uniform boundedness of $\m^\circ_n$ exactly describes
the usual \UMD-property \eqref{eq:umd} restricted to
Walsh-Paley-martingales.

Obviously, we have
\[
\m_n^{\circ\circ}(T) \leq \m_n^\circ(T) \leq \m_n(T).
\]
Surprisingly, there holds also an estimate in the reverse direction.
\begin{theorem*}
  $\m_n(T) \leq 3 \m_n^{\circ\circ}(T)$.
\end{theorem*}

\section{Proofs}
\label{sec:proofs}

For $(h,i)\in\D$, we denote by $\phi_h^{(i)}$ the transformation of
$[0,1)$ that interchanges the intervals
\[
  \Delta_{h+1}^{(4i-2)} \txt{and} \Delta_{h+1}^{(4i-1)}.
\]
More formally
\[
  \phi_h^{(i)}(t):=
  \begin{cases}
    t+\frac1{2^{h+1}} & \txt{for $t\in\Delta_{h+1}^{(4i-2)}$,} \\[3pt]
    t-\frac1{2^{h+1}} & \txt{for $t\in\Delta_{h+1}^{(4i-1)}$,} \\[3pt]
    t & \txt{otherwise.}
  \end{cases}
\]
It turns out that
\[
  \chi_k^{(j)}\circ \phi_h^{(i)} =
  \begin{cases}
    \chi_k^{(j)} &
    \txt{if $k<h$ or $k=h$, $j\not=i$,} \\[4pt]
    \displaystyle \frac{\chi_{h+1}^{(2i-1)}+\chi_{h+1}^{(2i)}}{\sqrt2} &
    \txt{if $k=h$ and $j=i$,} \\[12pt]
    \chi_k^{(j)} &
    \txt{if $k=h+1$ and $j\not=2i-1,2i$,} \\[4pt]
    \chi_k^{(j^*)} &
    \txt{if $k>h+1$,}
  \end{cases}
\]
where $(j^*)$ is a permutation of $(1,\dots,2^{k-1})$.
See \cite{wen96} for a proof.

The most important property for our purpose is that whenever
\[
\tsprod
f{\chi_{h+1}^{(2i-1)}} = \tsprod f{\chi_{h+1}^{(2i)}} = 0
\]
it follows that
\begin{equation} \label{eq:main_property}
  \tsprod{f\circ\phi_h^{(i)}}{\chi_h^{(i)}}=0 \txt{and} \tsprod
  {f\circ \phi_h^{(i)}}{\chi_{h+1}^{(2i-1)}} = \tsprod {f\circ
    \phi_h^{(i)}}{\chi_{h+1}^{(2i)}} = \frac{\tsprod
    f{\chi_h^{(i)}}}{\sqrt2}.
\end{equation}
In other words, the Haar-Fourier coefficient of a function $f$ with
respect to the index $(h,i)$ is shifted up one level and distributed to
the indices $(h+1,2i-1)$ and $(h+1,2i)$.

The basic idea of the proof is contained in the following proposition.
\begin{proposition}
  \label{prop:1}
  $\m_n(T) \leq \m_{2n}^{\circ\circ}(T)$.
\end{proposition}
\pf
For a $\D_1^n$-tuple $(x_k^{(j)})$ write
\[
f := \sum_{(k,j)\in\D_1^n} \!\! x_k^{(j)} \chi_k^{(j)}
\txt{and}
f^\epsilon := \sum_{(k,j)\in\D_1^n} \!\! \epsilon_k^{(j)} x_k^{(j)}
\chi_k^{(j)}.
\]
First, we want to find a transformation $\psi_1\colon[0,1)\to[0,1)$
such that the spectrum of $f\circ\psi_1$ is concentrated on the odd
levels, i.~e.
\[
\tsprod{f\circ\psi_1}{\chi_{2k}^{(j)}} = 0
\txt{for all $(2k,j)\in\D$.}
\]
Indeed, using the composition of all $\phi_n^{(j)}$ with
$j=1,\dots,2^{n-1}$, we shift the whole level $\D_n$ of the spectrum
of $f$ to the level $\D_{n+1}$.  Repeating this process of `spreading'
$\spec(f)$ successively on the levels $n+1,n+2,\dots,2n-2$ we move the
$n$-th level of $\spec(f)$ to the level $\D_{2n-1}$. In a similar
manner, we next move the $(n-1)$-st level to $\D_{2n-3}$ and so on.
So that finally
\[
\tsprod{f\circ\psi_1}{\chi_{2k}^{(j)}} = 0,
\]
as required.

Treating $f^\epsilon$ in the same way, we get that
\[
\tsprod{f^\epsilon\circ\psi_1}{\chi_{2k}^{(j)}} = 0
\]
and
\[
\tsprod{f^\epsilon\circ\psi_1}{\chi_{2k-1}^{(j)}}
= \delta_{2k-1}^{(j)}
\tsprod{f\circ\psi_1}{\chi_{2k-1}^{(j)}},
\]
where $\delta_{2k-1}^{(j)}=\pm1$ are signs that depend on the initial
signs $(\epsilon_k^{(j)})$ only.

We now construct a second transformation $\psi_2$ as composition of
all those transformations $\phi_{2k-1}^{(j)}$ for which
$\delta_{2k-1}^{(j)}=+1$. Since
\[
\tsprod{f\circ\psi_1}{\chi_{2k}^{(2j-1)}} =
\tsprod{f\circ\psi_1}{\chi_{2k}^{(2j)}} = 0,
\]
this moves all the plus signs to the even levels and leaves the minus
signs on the odd levels.  Letting $\psi:=\psi_2\circ\psi_1$, it
follows that
\[ \label{eq:pm1}
  \tsprod{f^\epsilon\circ\psi}{\chi_k^{(j)}}
  = (-1)^k
  \tsprod{f\circ\psi}{\chi_k^{(j)}}.
\]
Hence, the definition of $\m_{2n}^{\circ\circ}(T)$ yields
\[
\| Tf^\epsilon\circ\psi\|_2 \leq \m_{2n}^{\circ\circ}(T)
\|f\circ\psi\|_2.
\]
This completes the proof of Proposition \ref{prop:1}, since
\[
\| Tf^\epsilon\circ\psi\|_2 = \|Tf^\epsilon\|_2 =
\bigg\|
\sum_{(k,j)\in\D_1^n} \!\! \epsilon_k^{(j)} Tx_k^{(j)} \chi_k^{(j)}
\bigg\|_2
\]
and
\[
\| f\circ\psi\|_2 = \|f\|_2 =
\bigg\|
\sum_{(k,j)\in\D_1^n} \!\! x_k^{(j)} \chi_k^{(j)}
\bigg\|_2.\qed
\]

Next, we show that the sequence $\m_n^{\circ\circ}(T)$ behaves quite
regularly.
\begin{proposition}
  \label{prop:2}
  $\m_{2n}^{\circ\circ}(T) \leq 3 \m_n^{\circ\circ}(T)$.
\end{proposition}
\pf
Writing $\D_1^{2n}$ as the union of its lower part $\D_1^n$ and its
upper part $\D_{n+1}^{2n}$, we obtain
\[
\bigg\| \sum_{(k,j)\in\D_1^{2n}} \!\! (-1)^k Tx_k^{(j)}\chi_k^{(j)}
\bigg\|_2 \leq L + U,
\]
where
\[
L := \bigg\| \sum_{(k,j)\in\D_1^n} \!\! (-1)^k Tx_k^{(j)}\chi_k^{(j)}
\bigg\|_2 \txt{and}
U := \bigg\| \sum_{(k,j)\in\D_{n+1}^{2n}} \!\! (-1)^k Tx_k^{(j)}\chi_k^{(j)}
\bigg\|_2.
\]
Obviously
\[
L \leq \m_n^{\circ\circ}(T) \bigg\| \sum_{(k,j)\in\D_1^n}\!\!
x_k^{(j)}\chi_k^{(j)} \bigg\|_2,
\]
and by \eqref{eq:cond_exp}, we get
\[
\bigg\| \sum_{(k,j)\in\D_1^n} \!\! x_k^{(j)}\chi_k^{(j)} \bigg\|_2
\leq
\bigg\| \sum_{(k,j)\in\D_1^{2n}} \!\! x_k^{(j)}\chi_k^{(j)} \bigg\|_2.
\]
To estimate $U$, we use the `self-similarity' of the Haar functions.
Write $\D_{n+1}^{2n}$ as the disjoint union of its subtrees
\[
  \S_i := \{(k,j) \in \D_{n+1}^{2n} \colon
  j=(i-1)2^{k-n-1}+1,\dots,i2^{k-n-1}\}.
\]
Then the map
\[
  (k,j) \mapsto (k',j') := (k-n, j-(i-1)2^{k-n-1})
\]
defines a bijection of $\S_i$ and $\D_1^n$. Moreover, we have
\begin{equation} \label{eq:haar_properties}
  \chi_k^{(j)}(\tfrac{t+i-1}{2^n}) =
  \begin{cases}
    2^{n/2} \chi_{k'}^{(j')}(t) & \txt{if $(k,j)\in\S_i$,}
    \\[6pt]
    0 & \txt{otherwise.}
  \end{cases}
\end{equation}
Hence for
\begin{eqnarray*}
  U_i
  & := &
  \Bigg(
  \int_{\Delta_n^{(i)}}
  \bigg\| \sum_{(k,j)\in\D_{n+1}^{2n}} \!\! (-1)^k Tx_k^{(j)}\chi_k^{(j)}(t)
  \bigg\|^2 dt
  \Bigg)^{1/2} \\
  & = &
  \Bigg( \frac1{2^n}
  \int_0^1
  \bigg\| \sum_{(k,j)\in\S_i} \!\! (-1)^k
  Tx_k^{(j)}\chi_k^{(j)}(\tfrac{t+i-1}{2^n})
  \bigg\|^2 dt
  \Bigg)^{1/2} \\
  & = &
  \Bigg(
  \int_0^1
  \bigg\| \sum_{(k,j)\in\S_i} \!\! (-1)^k
  Tx_k^{(j)}\chi_{k'}^{(j')}(t)
  \bigg\|^2 dt
  \Bigg)^{1/2},
\end{eqnarray*}
we get
\begin{equation}
  U_i
  \leq
  \m_n^{\circ\circ}(T)
  \Bigg(
  \int_0^1
  \bigg\| \sum_{(k,j)\in\S_i}\!\!
  x_k^{(j)}\chi_{k'}^{(j')}(t)
  \bigg\|^2 dt
  \Bigg)^{1/2}.
  \label{eq:est_U.1}
\end{equation}
Using \eqref{eq:haar_properties} again, we obtain that
\begin{equation} \label{eq:est_U.2}
  \Bigg(
  \int_0^1
  \bigg\| \sum_{(k,j)\in\S_i}\!\!
  x_k^{(j)}\chi_{k'}^{(j')}(t)
  \bigg\|^2 dt
  \Bigg)^{1/2}
  =
  \Bigg(\!
  \int_{\Delta_n^{(i)}}
  \bigg\| \sum_{(k,j)\in\D_{n+1}^{2n}}
  \!\! x_k^{(j)}\chi_k^{(j)}(t)
  \bigg\|^2 dt
  \Bigg)^{1/2}\!.
\end{equation}
Putting \eqref{eq:est_U.1} and \eqref{eq:est_U.2} together yields
\[
U=\Big(\sum_{i=1}^{2^n} U_i^2 \Big)^{1/2} \leq
\m_n^{\circ\circ}(T)
\Bigg(
\int_0^1
\bigg\| \sum_{(k,j)\in\D_{n+1}^{2n}}
\!\! x_k^{(j)}\chi_k^{(j)}(t)
\bigg\|^2 dt
\Bigg)^{1/2}.
\]
Finally, again by \eqref{eq:cond_exp} we have
\[
\Bigg(
\int_0^1
\bigg\| \sum_{(k,j)\in\D_{n+1}^{2n}}
\!\! x_k^{(j)}\chi_k^{(j)}(t)
\bigg\|^2 dt
\Bigg)^{1/2} \leq 2
\Bigg(
\int_0^1
\bigg\| \sum_{(k,j)\in\D_1^{2n}}
\!\! x_k^{(j)}\chi_k^{(j)}(t)
\bigg\|^2 dt
\Bigg)^{1/2}.
\]
This completes the proof of Proposition \ref{prop:2}.\qed

The theorem is now an immediate consequence of Propositions
\ref{prop:1} and \ref{prop:2}.


\vspace{1ex}
\noindent
{\sc Mathematisches Institut, FSU Jena, 07740 Jena, Germany}\\
{\it E--mail:} {\tt wenzel@minet.uni-jena.de}

\end{document}